# On the Simes inequality and its generalization[*]

**Sanat K. Sarkar**[1]

*Temple University*

**Abstract:** The Simes inequality has received considerable attention recently because of its close connection to some important multiple hypothesis testing procedures. We revisit in this article an old result on this inequality to clarify and strengthen it and a recently proposed generalization of it to offer an alternative simpler proof.

## 1. Introduction

Let $X_{1:n} \leq \cdots \leq X_{n:n}$ denote the $n$ ordered values of a set of random variables $X_1, \ldots, X_n$. Assuming that the $X_i$'s are continuous with a common cdf $F$, Simes [36] considered the following inequality involving the ordered values $U_{1:n} \leq \cdots \leq U_{n:n}$ of $U(0,1)$ random variables $U_i = F(X_i)$, $i = 1, \ldots, n$,

(1.1) $$P\{U_{i:n} \geq i\alpha/n,\ i = 1, \ldots, n\} \geq 1 - \alpha,$$

while proposing a modification of the Bonferroni test for the intersection of null hypotheses. He conjectured that the inequality holds under positive dependence of the $X_i$'s (equivalently, the $U_i$'s), having proved the equality under the independence and numerically verified the inequality for some specific positively dependent multivariate distributions. Sarkar and Chang [29] proved this conjecture for a class of positively dependent multivariate distributions that are conditionally iid with a distribution that is stochastically increasing in the value of the conditioning variable, generalizing the work in Hochberg and Rom [10] and Samuel-Cahn [20] who first attempted to prove the conjecture in the bivariate case. Sarkar [22] later established the conjecture for a class of multivariate totally positive of order two (MTP$_2$) distributions that is larger than the one considered in Sarkar and Chang [29]. A careful study of Sarkar's [22] proof of the conjecture of course reveals that it holds for a slightly larger class of positively dependent multivariate distributions. This larger class is characterized by the following condition:

**Condition 1.1.** $E\{\phi(X_1, \ldots, X_n)|X_i\}$ is nondecreasing (or nonincreasing) in $X_i$ for each $i = 1, \ldots, n$, and any nondecreasing (or nonincreasing) function $\phi(X_1, \ldots, X_n)$.

---

[*]Supported by NSF Grant DMS-06-03868.
[1]Department of Statistics, Fox School of Business and Management, Temple University, Philadelphia, PA 19122, USA, e-mail: sanat@temple.edu

By a nondecreasing or nonincreasing function of more than one variable, in Condition 1.1 and elsewhere in the paper, we mean this to be so coordinatewise. The condition has been referred to as the positive dependence through stochastic ordering (PDS) condition by Block et al. [3].

Thus, Theorem 3.1 in Sarkar [22] establishing the Simes conjecture can be rephrased in its slightly improved form as follows, which we will refer to as the Simes inequality.

**Simes Inequality.** *Let $X_1, \ldots, X_n$ be a set of PDS continuous random variables with $F_i$ as the marginal cdf of $X_i$, $i = 1, \ldots, n$. Then, for any fixed $-\infty < a_1 \leq \cdots \leq a_n < \infty$,*

$$(1.2) \qquad P\{X_{1:n} \geq a_1, \ldots, X_{n:n} \geq a_n\} \geq 1 - \frac{1}{n}\sum_{i=1}^{n} F_i(a_n),$$

*if $j^{-1}F_i(a_j)$ is nondecreasing in $j = 1, \ldots, n$ for all $i = 1, \ldots, n$. The equality in (1.2) holds when $j^{-1}F_i(a_j)$ is constant in $j$ for each $i$ and the $X_i$'s are independent.*

Since the PDS property is invariant under co-monotone transformations of the $X_i$'s, the Simes inequality can be equivalently described with (1.2) replaced by the following:

$$(1.3) \qquad P\{X_{1:n} \leq b_1, \ldots, X_{n:n} \leq b_n\} \geq \frac{1}{n}\sum_{i=1}^{n} F_i(b_1),$$

for any fixed $-\infty < b_1 \leq \cdots \leq b_n < \infty$ such that $j^{-1}[1 - F_i(b_{n-j+1})]$ is nondecreasing in $j = 1, \ldots, n$ for all $i = 1, \ldots, n$. The equality in (1.3) holds when $j^{-1}[1 - F_i(b_{n-j+1})]$ is constant in $j$ for each $i$ and the $X_i$'s are independent.

The Simes inequality holds special importance in hypothesis testing. Besides theoretically validating the Type I error rate control of the Simes (1986) test, which has now been frequently used in place of the Bonferroni test in many scientific investigations (Dmitrienko et al. [7], Hommel et al. [11], Meng et al. [16], Neuhauser et al. [17], Rosenberg et al. [19], Somerville et al. [37] and Westfall and Krishen [38]), it provides theoretical basis for the familywise error rate (FWER) control under positive dependence of the commonly used Hochberg (1988) procedure for multiple testing, see, for example, Sarkar [22] and Sarkar and Chang [29]. Most importantly, it is closely linked to the inequality establishing the false discovery rate (FDR) control of the Benjamini and Hochberg [1] procedure (Benjamini and Yekutieli [2] and Sarkar [23, 24]), which as a FDR procedure has received the most attention so far in the multiple testing literature; see, for example, Sarkar [25, 26, 28] and Sarkar and Guo [30, 31] for references. The fact that the Simes inequality holds under Condition 1.1 is a by-product of the FDR control of the Benjamini-Hochberg procedure. Of course, the authors of the papers dealing with the Simes inequality were not aware of this condition being defined earlier as the PDS condition and referred to it as a special case of the positive regression dependence on subset (PRDS) condition under which the Benjamini-Hochberg procedure controls the FDR.

In this article, we will start with a generalized form of the Simes inequality that Sarkar [27] has recently obtained. First, we will provide an alternative simpler proof of this generalization. Then, we will go back to the original Simes inequality to clarify and strengthen an earlier result in Sarkar [22].

Sarkar [27] has generalized the Simes inequality by providing a lower bound for the probability $P\{X_{k:n} \geq a_k, \ldots, X_{n:n} \geq a_n\}$, for a fixed $1 \leq k \leq n$, in terms of the $k$th order joint distributions of the $X_i$'s in an attempt to generalize certain multiple testing procedures. He proved this generalization for MTP$_2$ distributions. We improve this work in this paper (Section 2) by offering an alternative simpler proof using a condition that is weaker than the MTP$_2$ condition.

To establish the Simes conjecture, that is, the original Simes inequality ($k = 1$), for both multivariate and absolute-valued multivariate central $t$ distributions when the associated multivariate normal has the same correlations, Sarkar [22] made use of a corollary (Corollary 3.1) attempting to extend his main theorem (Theorem 3.1) to certain scale mixtures of MTP$_2$ distributions. Unfortunately, while the main theorem is correct, there is a flaw in his proof of the corollary, as noted by Henry Block in a private communication. We will revisit this corollary in this article (Section 3) to clarify and at the same time strengthen it. More specifically, we will provide direct proofs of the Simes inequality for (i) multivariate $t$ distribution when the associated multivariate normal has nonnegative correlations and (ii) for absolute-valued multivariate $t$ distribution when the associated correlation matrix has a more general structure than just having equal correlations. These proofs will be based on ideas borrowed from Sarkar [23], although one can see Benjamini and Yekutieli [2] for a different proof, of course more complicated and given for a more general result.

## 2. Generalized Simes inequality

In this section, we give an alternative simpler proof of the generalized Simes inequality in Sarkar [27]. This simplification is achieved by offering a simpler proof of a supporting lemma on probability distribution of ordered random variables and by using the following condition for the distribution of **X** that is weaker, yet more directly applicable, than the MTP$_2$ condition. The following notation is used in the condition: $X_i^{-(i_1,\ldots,i_k)}$, $i = 1, \ldots, n-k$, are the components of the set $(X_1, \ldots, X_n) \setminus \{X_{i_1}, \ldots, X_{i_k}\}$.

**Condition 2.1.** For every $\{i_1, \ldots i_k\} \subseteq \{1, \ldots, n\}$ of size $k$ and any nondecreasing (or nonincreasing) function $\phi$, $E\{\phi(X_1^{-(i_1,\ldots,i_k)}, \ldots, X_n^{-(i_1,\ldots,i_k)}) | \max\{X_{i_1}, \ldots, X_{i_k}\} = z\}$ is nondecreasing (or nonincreasing) in $z$.

**Theorem 2.1** (Generalized Simes inequality). *Under Condition 2.1, we have*

$$(2.1) \quad P\{X_{k:n} \geq a_k, \ldots, X_{n:n} \geq a_n\} \geq 1 - \binom{n}{k}^{-1} \sum_{1 \leq i_1 < \cdots < i_k \leq n} F_{i_1 \ldots i_k}(a_n),$$

*with $F_{i_1 \ldots i_k}(x) = P\{\max_{1 \leq j \leq k} X_{i_j} \leq x\}$, for any fixed $1 \leq k \leq n$, where $a_k \leq \cdots \leq a_n$ are such that $\binom{j}{k}^{-1} F_{i_1 \ldots i_k}(a_j)$ is nondecreasing in $j = k, \ldots, n$ for every $\{i_1, \ldots, i_k\} \subseteq \{1, \ldots, n\}$ of cardinality $k$. The equality in (2.1) holds under the independence and when $\binom{j}{k}^{-1} F_{i_1 \ldots i_k}(a_j)$ is constant in $j = k, \ldots, n$ for every $\{i_1, \ldots, i_k\} \subseteq \{1, \ldots, n\}$ of cardinality $k$.*

The theorem will be proved using the following lemma. As mentioned above, this lemma, although proved before in Sarkar [27], will be given an alternative simpler proof here.

**Lemma 2.1.** *Given an increasing set of constants $a_1 \leq \cdots \leq a_n$, let*

$$(2.2) \qquad R_n = \max_{1 \leq i \leq n} \{i : X_{i:n} \leq a_i\}.$$

*Then,*

$$(2.3) \quad \begin{aligned} P(R_n \geq k) &= \sum_{1 \leq i_1 < \cdots < i_k \leq n} P\left\{\max_{1 \leq j \leq k} X_{i_j} \leq a_k\right\} - \sum_{r=k+1}^{n} \sum_{1 \leq i_1 < \cdots < i_k \leq n} \\ & E\left(P\left\{R_{n-k}^{-(i_1,\ldots,i_k)} \geq r-k \mid X_{i_1}, \ldots, X_{i_k}\right\} \left[\binom{r-1}{k}^{-1} I\left\{\max_{1 \leq j \leq k} X_{i_j} \leq a_{r-1}\right\} \right. \right. \\ & \left. \left. - \binom{r}{k}^{-1} I\left\{\max_{1 \leq j \leq k} X_{i_j} \leq a_r\right\}\right]\right), \end{aligned}$$

*where*

$$(2.4) \qquad R_{n-k}^{-(i_1,\ldots,i_k)} = \max_{1 \leq i \leq n-k} \left\{i : X_{i:n-k}^{-(i_1,\ldots,i_k)} \leq a_{k+i}\right\},$$

*with $X_{i:n-k}^{-(i_1,\ldots,i_k)}$, $i = 1, \ldots, n-k$, being the $n-k$ ordered values of the set $(X_1, \ldots, X_n) \setminus \{X_{i_1}, \ldots, X_{i_k}\}$.*

*Proof.* Given that $R_n = r$, where $k \leq r \leq n$, $R_n$ can be expressed as $R_n = \sum_{i=1}^n I(X_i \leq a_r)$. Hence, we have

$$(2.5) \quad \begin{aligned} R_n(R_n - 1) \cdots (R_n - k + 1) &= \sum_{1 \leq i_1 \neq \cdots \neq i_k \leq n} I(\max_{1 \leq j \leq k} X_{i_j} \leq a_r) \\ &= k! \sum_{1 \leq i_1 < \cdots < i_k \leq n} I(\max_{1 \leq j \leq k} X_{i_j} \leq a_r). \end{aligned}$$

Therefore, for $k \leq r \leq n$, we have

$$(2.6) \quad \begin{aligned} I(R_n = r) &= \frac{R_n(R_n - 1) \cdots (R_n - k + 1)}{r(r-1) \cdots (r-k+1)} I(R_n = r) \\ &= \binom{r}{k}^{-1} \sum_{1 \leq i_1 < \cdots < i_k \leq n} I(\max_{1 \leq j \leq k} X_{i_j} \leq a_r, R_n = r), \end{aligned}$$

which yields

$$(2.7) \quad \begin{aligned} I(R_n \geq k) &= \sum_{r=k}^{n} I(R_n = r) \\ &= \sum_{r=k}^{n} \binom{r}{k}^{-1} \sum_{1 \leq i_1 \leq \cdots \leq i_k \leq n} I\left(\max_{1 \leq j \leq k} X_{i_j} \leq a_r, R_n = r\right) \\ &= \sum_{r=k}^{n} \binom{r}{k}^{-1} \sum_{1 \leq i_1 \leq \cdots \leq i_k \leq n} I\left(\max_{1 \leq j \leq k} X_{i_j} \leq a_r, R_{n-k}^{-(i_1,\ldots,i_k)} = r - k\right). \end{aligned}$$

Suppressing the superscript in $R_{n-k}$ and using $I(R_{n-k} = r - k) = I(R_{n-k} \geq r - k) - I(R_{n-k} \geq r - k + 1)$ for $r = k, \ldots, n$ in (2.7), we get

$$
\begin{aligned}
(2.8) \quad & I(R_n \geq k) \\
&= \sum_{1 \leq i_1 \leq \cdots \leq i_k \leq n} I\left(\max_{1 \leq j \leq k} X_{i_j} \leq a_k\right) - \sum_{1 \leq i_1 \leq \cdots \leq i_k \leq n} \sum_{r=k+1}^{n} I(R_{n-k} \geq r - k) \\
& \left[\binom{r-1}{k}^{-1} I\left(\max_{1 \leq j \leq k} X_{i_j} \leq a_{r-1}\right) - \binom{r}{k}^{-1} I\left(\max_{1 \leq j \leq k} X_{i_j} \leq a_r\right)\right].
\end{aligned}
$$

Taking expectations of both sides in (2.8), we get the lemma. □

**Remark 2.1.** Before we proceed to prove Theorem 2.1, it is important to note that formula (2.3) can alternatively be written as

$$
\begin{aligned}
(2.9) \quad & P(R_n \geq k) \\
&= \binom{n}{k}^{-1} \sum_{1 \leq i_1 < \cdots < i_k \leq n} P\left\{\max_{1 \leq j \leq k} X_{i_j} \leq a_n\right\} - \sum_{r=k+1}^{n} \sum_{1 \leq i_1 < \cdots < i_k \leq n} \\
& E\left(P\left\{0 \leq R_{n-k}^{-(i_1,\ldots,i_k)} < r - k \mid X_{i_1}, \ldots, X_{i_k}\right\} \left[\binom{r}{k}^{-1} I\left\{\max_{1 \leq j \leq k} X_{i_j} \leq a_r\right\} \right.\right. \\
& \left.\left. - \binom{r-1}{k}^{-1} I\left\{\max_{1 \leq j \leq k} X_{i_j} \leq a_{r-1}\right\}\right]\right).
\end{aligned}
$$

*Proof of Theorem 2.1.* Since

$$
(2.10) \quad P\{X_{k:n} \geq a_k, \ldots, X_{n:n} \geq a_n\} = 1 - P\{R_n \geq k\},
$$

the theorem follows from the fact that, for every fixed $r = k + 1, \ldots, n$ and $\{i_1, \ldots, i_k\} \subseteq \{1, \ldots, n\}$, the expectation under the multiple summation in the right-hand side of (2.9) is greater than or equal to 0, which can be proved as follows. Define

$$
\psi_{r,k}(Z) = P\left\{0 \leq R_{n-k}^{-(i_1,\ldots,i_k)} < r - k \mid \max_{1 \leq j \leq k} X_{i_j} = Z\right\}.
$$

Then, the above expectation is

$$
\begin{aligned}
& E\left\{\psi_{r,k}(Z)\left[\binom{r}{k}^{-1} I(Z \leq a_r) - \binom{r}{k-1}^{-1} I(Z \leq a_{r-1})\right]\right\} \\
=& E\left\{\psi_{r,k}(Z) I(Z \leq a_r)\left[\binom{r}{k}^{-1} - \binom{r}{k-1}^{-1} I(Z \leq a_{r-1})\right]\right\} \\
\geq& \frac{E\{\psi_{r,k}(Z) I(Z \leq a_r)\}}{P(Z \leq a_r)} \left[\binom{r}{k}^{-1} P(Z \leq a_r) - \binom{r}{k-1}^{-1} P(Z \leq a_{r-1})\right] \\
\geq& 0,
\end{aligned}
$$

as $\psi_{r,k}(Z)$ is a nondecreasing function of $Z$, because of Condition 2.1 and $\{0 \leq R_{n-k}^{-(i_1,\ldots,i_k)} < r - k\}$ being a nondecreasing set, and $\binom{r}{k}^{-1} - \binom{r}{k-1}^{-1} I(Z \leq a_{r-1})$ is also a nondecreasing function of $Z$. □

**Remark 2.2.** One can get an equivalent statement of Theorem 2.1, generalizing (1.3), the equivalent version of the Simes inequality (1.2), by replacing each $X_i$ by, say $-X_i$. But, since unlike the PDS condition, Condition 2.1 is not invariant under co-monotone transformations, one needs to have a condition different from Condition 2.1 in this alternative statement of Theorem 2.1. It is not difficult to see the following alternative versions of Condition 2.1 and Theorem 2.1:

**Condition 2.1\*.** For every $\{i_1, \ldots i_k\} \subseteq \{1, \ldots, n\}$ of size $k$ and any nondecreasing (or nonincreasing) function $\phi$, $E\{\phi(X_1^{-(i_1,\ldots,i_k)}, \ldots, X_n^{-(i_1,\ldots,i_k)}) | \min\{X_{i_1}, \ldots, X_{i_k}\} = z\}$ is nondecreasing (or nonincreasing) in $z$.

**Theorem 2.1\*.** *Under Condition 2.1\*, we have*

$$P\{X_{1:n} \leq b_1, \ldots, X_{n-k+1:n} \leq b_{n-k+1}\}$$
$$(2.11) \qquad \geq \binom{n}{k}^{-1} \sum_{1 \leq i_1 < \cdots < i_k \leq n} G_{i_1 \ldots i_k}(b_1),$$

*with $G_{i_1\ldots i_k}(x) = P\{\min_{1 \leq j \leq k} X_{i_j} \leq x\}$, for any fixed $1 \leq k \leq n$, where $b_1 \leq \cdots \leq b_{n-k+1}$ are such that $\binom{j}{k}^{-1}[1 - G_{i_1\ldots i_k}(b_{n-j+1})]$ is nondecreasing in $j = k, \ldots, n$ for every $\{i_1, \ldots, i_k\} \subseteq \{1, \ldots, n\}$ of cardinality $k$. The equality in (2.11) holds under the independence and when $\binom{j}{k}^{-1}[1 - G_{i_1\ldots i_k}(b_{n-j+1})]$ is constant in $j = k, \ldots, n$ for every $\{i_1, \ldots, i_k\} \subseteq \{1, \ldots, n\}$ of cardinality $k$.*

**Remark 2.3.** Condition 2.1 or 2.1\*, for $1 < k \leq n$, is more restrictive than the PDS condition. So, Theorem 2.1 or 2.1\* holds for a smaller class of distributions than the one for which the Simes original inequality holds. In fact, we have the following lemma providing a class of distributions for which the generalized Simes inequality holds. This will be proved using properties of $\text{MTP}_2$ distributions discussed in Karlin and Rinott [12].

**Lemma 2.2.** *Both Conditions 2.1 and 2.1\* are satisfied when $(X_1, \ldots, X_n)$ has a symmetric $\text{MTP}_2$ distribution.*

*Proof.* Let $(X_1, \ldots, X_n) \sim f(x_1, \ldots, x_n)$, which is symmetric and $\text{MTP}_2$, and $Y = \max\{X_1, \ldots, X_k\}$. The joint density of $Y, X_{k+1}, \ldots, X_n$ is given by

$$g(y, x_{k+1}, \ldots, x_n)$$
$$(2.12) \quad = k \int \cdots \int f(y, x_2, \ldots, x_k, x_{k+1}, \ldots, x_n) \prod_{i=2}^{k} I(x_i \leq y) \prod_{i=2}^{k} dx_i.$$

This is $\text{MTP}_2$, since $f$ and $\prod_{i=2}^{k} I(x_i \leq y)$ are both $\text{MTP}_2$, their product is $\text{MTP}_2$, and so is the above integral. Therefore, $(Y, X_1, \ldots, X_n)$ satisfies the positive regression dependence condition, which implies that

$$E\left[\phi(X_{k+1}, \ldots, X_n) \big| \max\{X_1, \ldots, X_k\} = z\right]$$

is nondecreasing (nonincreasing) in $z$ for any nondecreasing (nonincreasing) function $\phi$, which is Condition 2.1 in this case.

With $Y = \min\{X_1, \ldots, X_k\}$, $(Y, X_{k+1}, \ldots, X_n)$ is also jointly $\text{MTP}_2$, which can be proved as before by replacing $\prod_{i=2}^{k} I(x_i \leq y)$ by $\prod_{i=2}^{k} I(x_i \geq y)$ in (2.12) and using the fact that $\prod_{i=2}^{k} I(x_i \geq y)$ is also $\text{MTP}_2$. So, Condition 2.1\* is also satisfied. □

Lemma 2.2 now yields the following corollary to both Theorems 2.1 and 2.1*.

**Corollary 2.1.** *Suppose that $(X_1, \ldots, X_n)$ have a symmetric $MTP_2$ distribution.*

*(a) Let $F_k(x)$ be the common cdf of the maximum of any $k$ of the $n$ $X_i$'s and $a_k \leq \cdots \leq a_n$ be such that $\binom{j}{k}^{-1} F_k(a_j)$ is nondecreasing in $j = k, \ldots, n$. Then,*

$$P\{X_{k:n} \geq a_k, \ldots, X_{n:n} \geq a_n\} \geq 1 - F_k(a_n). \tag{2.13}$$

*The equality holds under the independence and when $\binom{j}{k}^{-1} F_k(a_j)$ is constant in $j = k, \ldots, n$.*

*(b) Let $G_k(x)$ be the common cdf of the minimum of any $k$ of the $n$ $X_i$'s and $b_1 \leq \cdots \leq b_{n-k+1}$ be such that $\binom{j}{k}^{-1}[1 - G_k(b_{n-j+1})]$ is nondecreasing in $j = k, \ldots, n$. Then,*

$$P\{X_{1:n} \leq b_1, \ldots, X_{n-k+1:n} \leq b_{n-k+1}\} \geq G_k(b_1), \tag{2.14}$$

*The equality holds under the independence and when $\binom{j}{k}^{-1}[1 - G_k(b_{n-j+1})]$ is constant in $j = k, \ldots, n$.*

*Proof.* The expression (2.10) in this case reduces to

$$P\{X_{k:n} \geq a_k, \ldots, X_{n:n} \geq a_n\}$$
$$= 1 - F_k(a_n) + \binom{n}{k} \sum_{r=k+1}^{n}$$
$$E\left\{\psi_{r,k}(Y) \left[\binom{r}{k}^{-1} I(Y \leq a_r) - \binom{r-1}{k}^{-1} I(Y \leq a_{r-1})\right]\right\}, \tag{2.15}$$

where

$$\psi_{r,k}(Y) = E\{\phi_{r,k}(X_{k+1}, \ldots, X_n)|Y\},$$

with $\phi_{r,k}(X_{k+1}, \ldots, X_n)$ as the indicator function of the event

$$\left\{0 \leq R_{n-k}^{-(1,\ldots,k)} < r - k\right\}$$
$$\equiv \left\{X_{r-k:n-k}^{-(1,\ldots,k)} \geq a_r, \ldots, X_{n-k:n-k}^{-(1,\ldots,k)} \geq a_n\right\}. \tag{2.16}$$

Since Condition 2.1 is now satisfied (because of Lemma 2.2) and the indicator function $\phi_{r,k}$ is a nondecreasing function of $(X_{k+1}, \ldots, X_n)$, $\psi_{r,k}(Y)$ is a nondecreasing function of $Y$, which proves (2.13) as in the proof of Theorem 2.1.

The part (b) of the corollary can be similarly proved by noting that

$$P\{X_{1:n} \leq b_1, \ldots, X_{n-k+1:n} \leq b_{n-k+1}\}$$
$$= G_k(b_1) + \binom{n}{k} \sum_{r=k+1}^{n}$$
$$E\left\{\psi_{r,k}(Y) \left[\binom{r}{k}^{-1} I(Y \geq b_{n-r+1}) - \binom{r-1}{k}^{-1} I(Y \geq b_{n-r+2})\right]\right\}, \tag{2.17}$$

where now $Y = \min\{X_1, \ldots, X_k\}$ and

$$\psi_{r,k}(Y) = E\{\phi_r(X_{k+1}, \ldots, X_n)|Y\},$$

for some nonincreasing function $\phi_{r,k}$. The inequality (2.14) follows because $\psi_{r,k}(Y)$ is nonincreasing in $Y$. □

The covariance matrix $\boldsymbol{\Sigma}$ of a symmetric multivariate normal distribution with a common non-negative correlation satisfies each of the properties: (i) the off-diagonals of $-\boldsymbol{\Sigma}^{-1}$ are non-negative and (ii) the off-diagonals of $-D\boldsymbol{\Sigma}^{-1}D$ are non-negative for some diagonal matrix $D$ with diagonal entries $\pm 1$, which are the conditions, respectively, for multivariate normal, $N_n(\boldsymbol{\mu}, \boldsymbol{\Sigma})$, and absolute-valued zero-mean multivariate normal, $|N_n(\mathbf{0}, \boldsymbol{\Sigma})|$, to be MTP$_2$; see, for example, Karlin and Rinott [12, 13]. Thus, we have the following:

**Proposition 2.1.** *The generalized Simes inequality holds for both symmetric and absolute-valued zero-mean symmetric multivariate normal distributions with a common non-negative correlation.*

If the above distributions are studentized based on an independent chi-square random variable, the resulting multivariate and absolute-valued multivariate $t$ distributions may not retain the MTP$_2$ property. So, it becomes unclear if the generalized Simes' inequality still holds for these distributions, although the Simes original inequality does. In fact, as we will show in the next section, the associated covariance matrices for these $t$ distributions do not have to be symmetric for the original Simes inequality to hold.

## 3. Simes inequalities for $t$ distributions

In this section, we will revisit the Simes inequalities for multivariate and absolute-valued multivariate $t$ distributions to clarify and strengthen previous related work. To be more specific, we have the following theorem.

**Theorem 3.1.** *Let $T_i = Z^{-1}X_i$, $i = 1, \ldots, n$, where $(X_1, \ldots, X_n) \sim N_n(\mathbf{0}, \boldsymbol{\Sigma})$ with the diagonal entries of $\boldsymbol{\Sigma}$ being 1, and is independent of $Z \sim \chi_\nu/\sqrt{\nu}$. Then, the Simes inequality (1.2) [or (1.3)] holds (i) for the $T_i$'s if $a_n \leq 0$ (or $b_1 \geq 0$) and the off-diagonals of $\boldsymbol{\Sigma}$ are non-negative and (ii) for the $|T_i|$'s if the off-diagonals of $-D\boldsymbol{\Sigma}^{-1}D$ are non-negative for some diagonal matrix $D$ with diagonal entries $\pm 1$.*

Before we proceed to prove this theorem, it is important to re-emphasize that it is the PDS condition that drives the Simes inequality, and hence there are distributions, not necessarily MTP$_2$, for which the inequality holds. A case in point is multivariate normal with nonnegative correlations. Its PDS property follows easily from the fact that the conditional means given any $X_i$ are increasing in that $X_i$, even though it is not MTP$_2$ unless the off-diagonals of $-\boldsymbol{\Sigma}^{-1}$ are also nonnegative [12]. For absolute-valued multivariate normal distribution, of course, the PDS property does not follow that easily unless the MTP$_2$ property is invoked, and that holds when the off-diagonals of $-D\boldsymbol{\Sigma}^{-1}D$ are non-negative for some diagonal matrix $D$ with diagonal entries $\pm 1$ [12, 13].

Having proved the Simes inequality for MTP$_2$ distributions, Sarkar [22] attempted to prove it (in Corollary 3.1) for scale-mixtures of certain symmetric MTP$_2$ distributions before discussing that the inequality holds for symmetric multivariate $t$ and absolute-valued symmetric multivariate $t$ distributions. Unfortunately, as noted in the introduction, there is a flaw in his proof of the corollary. Nevertheless, while the truth of the corollary becomes an open issue at this point, it is important to emphasize that the Simes inequalities for these $t$ distributions that the corollary intends to prove still hold, as noted when dealing with similar inequalities arising in the context of the FDR control of the Benjamini-Hochberg procedure (Benjamini and Yekutieli [2] and Sarkar [23]). In Benjamini and Yekutieli [2], a general result on the PDS property of certain scale-mixtures of PDS distributions is given, from which one can see that the Simes inequality can be extended from a multivariate or absolute-valued multivariate normal to the corresponding multivariate or absolute-valued multivariate $t$ distribution. While this result is important in its own right, its proof, however, seems complicated. In fact, one can avoid it while extending the Simes inequality from multivariate normal to the corresponding multivariate $t$ distribution and, instead, apply an independence result of normal distribution. This has been briefly pointed out in Sarkar [23], of course, in the context of the FDR. We elaborate this point here in the context of Simes inequality, thereby clarifying and strengthening the inequalities for symmetric multivariate $t$ and absolute-valued symmetric multivariate $t$ distributions discussed in Sarkar [22]. We are now giving alternative and direct proofs of these inequalities with covariance matrices that are not necessarily symmetric.

*Proof of Theorem 3.1.* From (2.10), we have

$$
\begin{aligned}
& P\{T_{1:n} \geq a_1, \ldots, T_{n:n} \geq a_n\} \\
(3.1) \quad &= 1 - F(a_n) + \sum_{i=1}^{n} \sum_{r=2}^{n} E\left\{\psi_r(T_i)\left[\frac{I(T_i \leq a_r)}{r} - \frac{I(T_i \leq a_{r-1})}{r-1}\right]\right\},
\end{aligned}
$$

where $F$ is the common cdf of each $T_i$ and

$$
\psi_r(T_i) = P\left\{T^{-(i)}_{r-1:n-1} \geq a_r, \ldots, T^{-(i)}_{n-1:n-1} \geq a_n \mid T_i\right\}.
$$

We will prove in the following that each expectation under the double summation in (3.1) is greater than or equal to zero if $F(a_i)/i$ is nondecreasing in $i$ as long as $a_n \leq 0$, which will prove the Simes inequality (1.2).

Let us consider the expectation for $i = n$, and assume without any loss of generality that $Z \sim \chi_\nu$. Then, this expectation is given by

$$
(3.2) \quad E\left(P\left\{X^{-(n)}_{r-1:n-1} \geq a_r Z, \ldots, X^{-(n)}_{n-1:n-1} \geq a_n Z \mid X_n, Z\right\} \left[\frac{I(X_n \leq a_r Z)}{r} - \frac{I(X_n \leq a_{r-1} Z)}{r-1}\right]\right).
$$

Let

$$
(3.3) \quad g(x,z) = P\left\{X^{-(n)}_{r-1:n-1} \geq a_r z, \ldots, X^{-(n)}_{n-1:n-1} \geq a_n z \mid X_n = x\right\},
$$

where $z > 0$. Then, the expectation in (3.2) can be rewritten in terms of independent random variables $Z_n^* = \sqrt{Z^2 + X_n^2}$ and $T_n^* = T_n/\sqrt{1+T_n^2}$ as follows

$$
\begin{aligned}
& E\left\{g\left(Z_n^* T_n^*, Z_n^*\sqrt{1-T_n^{*2}}\right)\left[\frac{I(T_n^* \leq a_r^*)}{r} - \frac{I(T_n^* \leq a_{r-1}^*)}{r-1}\right]\right\} \\
(3.4) \quad &= E\left\{h(T_n^*)\left[\frac{I(T_n^* \leq a_r^*)}{r} - \frac{I(T_n^* \leq a_{r-1}^*)}{r-1}\right]\right\},
\end{aligned}
$$

where $a_r^* = a_r/\sqrt{1+a_r^2}$ and

$$
(3.5) \quad h(T_n^*) = E\left\{g\left(Z_n^* T_n^*, Z_n^*\sqrt{1-T_n^{*2}}\right) \mid T_n^*\right\}.
$$

Since the $X_i$'s are PDS and $a_i$'s are assumed negative, the probability in (3.3) is nondecreasing in $(x, z)$, implying that (3.5) is a nondecreasing function of $T_n^*$ as long as $T_n^* < 0$. Hence, the expectation in (3.4) is greater than or equal to

$$\frac{E\{h(T_n^*)I(T_n^* \leq a_r^*)\}}{P\{T_n^* \leq a_r^*\}} \left[\frac{P\{T_n^* \leq a_r^*\}}{r} - \frac{P\{T_n^* \leq a_{r-1}^*\}}{r-1}\right]$$

$$(3.6) \quad = \frac{E\{h(T_n^*)I(T_n^* \leq a_r^*)\}}{P\{T \leq a_r\}} \left[\frac{P\{T_n \leq a_r\}}{r} - \frac{P\{T_n \leq a_{r-1}\}}{r-1}\right],$$

which is greater than or equal to zero. Thus, the Simes inequality (1.2) holds. The version (1.3) of the Simes inequality can be similarly proved with positive $b_i$'s.

To prove Theorem 3.1(ii), we continue with the same arguments as above replacing $X_i$ (or $T_i$) by $|X_i|$ (or $|T_i|$). The $|X_i|$'s are PDS because of being $\mathrm{MTP}_2$ under the assumed condition on the covariance matrix. So, the function $g(x, z)$ now is nondecreasing in $x$ and is nonincreasing in $z$, implying that the function $h(|T_n^*|)$ continues to be an increasing function of $|T_n^*|$. The rest of the arguments remains same, completing the proof. □

## 4. Concluding remarks

The results discussed in this article basically are probability inequalities for the ordered components of a certain type of positively dependent random variables. More specifically, they provide bounds for joint probabilities of the ordered components of a set of random variables in terms of lower dimensional marginal distributions under a form of positive dependence among the variables. Our primary focus in this paper has been on these inequalities, rather than on discussing about the related Simes tests validated by these inequalities and their use in multiple testing procedures. Readers can see Cai and Sarkar [4, 5], Hochberg and Liberman [9], Krummenauer and Hommel [15], Rødland [18], Samuel-Cahn [21], Sen [32], Sen and Silvapulle [33], Seneta and Chen [34], Silvapulle and Sen [35] for the Simes test and Sarkar [27] for its generalization, in addition to those cited before.

Given two independent random samples $(X_1, \ldots, X_n)$ and $(Y_1, \ldots, Y_n)$ from two continuous populations $F$ and $G$ respectively, Lemma 2.1 has a potential application in developing a nonparametric test for testing the null hypothesis that $F$ and $G$ are equal versus $G$ is stochastically larger than $F$. Specifically, one can consider the statistic

$$T_n = \max_{1 \leq i \leq n} \{X_{i:n} \leq Y_{i:n}\},$$

where $X_{1:n} \leq \cdots \leq X_{n:n}$ and $Y_{1:n} \leq \cdots \leq Y_{n:n}$ are the order statistics corresponding to these samples, the probability distribution of which under the null hypothesis can be explicitly obtained using this lemma. A class of tests based on U-statistics with kernels based on sub-sample maximas is proposed in Deshpande and Kochar [6] and Kochar [14]. Perhaps one can propose and study tests based on a combination of some or all members of this class to increase the efficiency. The results obtained in this paper may be useful in finding the p-values of such tests.

theorem. *J. Statist. Plann. Inference* **82** 147–149. MR1736439
[22] SARKAR, S. K. (1998). Some probability inequalities for ordered $MTP_2$ random variables: a proof of the Simes conjecture. *Ann. Statist.* **26** 494–504. MR1626047
[23] SARKAR, S. K. (2002). Some results on false discovery rate in stepwise multiple testing procedures. *Ann. Statist.* **30** 239–257. MR1892663
[24] SARKAR, S. K. (2004). FDR-controlling stepwise procedures and their false negatives rates. *J. Statist. Plann. Inference* **125** 119–137. MR2086892
[25] SARKAR, S. K. (2006). False discovery and false nondiscovery rates in single-step multiple testing procedures. *Ann. Statist.* **34** 394–415. MR2275247
[26] SARKAR, S. K. (2007a). Stepup procedures controlling generalized FWER and generalized FDR. *Ann. Statist.* To appear. MR2382652
[27] SARKAR, S. K. (2007b). Generalizing Simes' test and Hochberg's stepup procedure. *Ann. Statist.* To appear. MR2387974
[28] SARKAR, S. K. (2007c). Two-stage stepup procedures controlling FDR. *J. Statist. Plann. Inference.* To appear. MR2384506
[29] SARKAR, S. K. AND CHANG, C.-K. (1997). The Simes method for multiple hypothesis testing with positively dependent test statistics. *J. Amer. Statist. Assoc.* **92** 1601–1608. MR1615269
[30] SARKAR, S. K. AND GUO, W. (2006). Procedures controlling generalized false discovery rate. Technical report, Temple Univ. Available at http://astro.temple.edu/~sanat/reports/GeneralizedFDR.pdf.
[31] SARKAR, S. K. AND GUO, W. (2007). On generalized false discovery rate. Unpublished manuscript.
[32] SEN, P. K. (1999). Some remarks on Simes-type multiple tests of significance. *J. Statist. Plann. Inference* **82** 139–145. MR1736438
[33] SEN, P. K. AND SILVAPULLE, M. J. (2002). An appraisal of some aspects of statistical inference under inequality constraints. *J. Statist. Plann. Inference* **107** 3–43. MR1927753
[34] SENETA, E. AND CHEN, J. (2005). Simple stepwise tests of hypotheses and multiple comparisons. *Internat. Statist. Rev.* **73** 21–34.
[35] SILVAPULLE, M. J. AND SEN, P. K. (2004). *Constrained Statistical Inference.* Wiley, New York. MR2099529
[36] SIMES, R. J. (1986). An improved Bonferroni procedure for multiple tests of significance. *Biometrika* **73** 751–754. MR0897872
[37] SOMERVILLE, M., WILSON, T., KOCH, G. AND WESTFALL, P. (2005). Evaluation of a weighted multiple comparison procedure. *Pharm. Statist.* **4** 7–13.
[38] WESTFALL, P. H. AND KRISHEN, A. (2001). Optimally weighted, fixed sequence and gatekeeper multiple testing procedures. *J. Statist. Plann. Inference* **99** 25–40. MR1858708